\documentstyle[12pt]{article}
\begin{document}
\author{
{\bf V. P. Akulov, }\\
{\it Theory Division, Kharkov Institute of Physics and Technology,}\\
{\it Kharkov 310108, Ukraine,}\\
{\bf V. V. Chitov}\\
{\it Department of Physics and Technology, Kharkov State University,}\\
{\it  Kharkov 310077, Ukraine}\\
{\bf Steven Duplij} \thanks{%
Alexander von Humboldt Fellow. On leave of absence from Theory Division,
Nuclear Physics Laboratory, Kharkov State University, Kharkov 310077,
Ukraine. Electronic address: duplij@physik.uni-kl.de. Internet address:
http://gluon.physik.uni-kl.de/\~{}duplij } \\
%EndAName
{\it Physics Department, University of Kaiserslautern},\\
{\it D-67653 Kaiserslautern, Germany}}
\title{{\bf DIFFERENTIAL CALCULUS ON
$q$-DEFORMED LIGHT-CONE
}}
\date{January 12, 1997}
\maketitle
\begin{abstract}
We propose the ``short'' version of $q$-deformed
differential calculus on the
light-cone using twistor representation.
The commutation relations between
coordinates and momenta are obtained. The quasi-classical limit introduced
gives an exact shape of the off-shell shifting.
\end{abstract}


\newpage
\section{Introduction}

The deformed differential calculus \cite{wes/zum3,wor} plays an important
role in physical applications \cite{are/vol1,mue-hoi}. In this connection
the light-cone approach attracts the special interest, because it is
employed for the description of massless particles, null-strings and
null-membranes \cite{gre/sch/wit}.The usual formalism of the bicovariant
differential calculi \cite{fad/pya,are/aru/med} and right (left) invariant
differential calculi \cite{aku/ger} are not applicable on the light-cone due
to vanishing of the corresponding quantum determinant. That is the reason
why the Maurer-Cartan forms cannot be defined here at all, and therefore it
is not possible to build the deformed differential calculus on the
light-cone in the standard way.

From general mathematical viewpoint some objects under consideration
belong to the class of
quantum semigroups \cite{dem} which consist of the standard quantum groups
and ideals (for abstract semigroup theory see \cite{cli/pre1} and for
application to supersymmetry see \cite{dup6,dup14}).

Here we construct a version of $q$-deformed differential calculus on the
light-cone which is mostly close to the standard one \cite{zup1,vak/kor}. We
therefore hope that the obtained calculus can be directly used in physical
applications.

\section{``Short'' differential calculus}

Let us remind the deformed differential calculi on the quantum plane $%
E_q\left( 2\right) $ \cite{wes/zum3} 
\begin{equation}
xy=qyx.  \label{q}
\end{equation}
In the standard case we have two solutions for $E_q\left( 2\right) $

a) 
\begin{equation}
\left\{ 
\begin{array}{lll}
\delta xx & = & q^2x\delta x, \\ 
\delta xy & = & qy\delta x+\left( q^2-1\right) x\delta y, \\ 
\delta yx & = & qx\delta y, \\ 
\delta yy & = & q^{-2}y\delta y;
\end{array}
\right.  \label{a}
\end{equation}

b)

\begin{equation}
\left\{ 
\begin{array}{lll}
\delta xx & = & q^2x\delta x, \\ 
\delta xy & = & q^{-1}y\delta x, \\ 
\delta yx & = & q^{-1}x\delta y+\left( q^{-2}-1\right) y\delta x, \\ 
\delta yy & = & q^{-2}y\delta y.
\end{array}
\right.  \label{b}
\end{equation}

It is well-known \cite{fad/res/tak} that there exists the automorphism of
the quantum plane (\ref{q}) 
\begin{equation}
\left( 
\begin{array}{c}
x\rightarrow y \\ 
y\rightarrow x \\ 
q\rightarrow q^{-1}
\end{array}
\right)  \label{auto}
\end{equation}
which is stipulated by the following property of the $R$-matrix 
\begin{equation}
R_q\rightarrow R_{q^{-1}}^{-1}.  \label{r}
\end{equation}

If we introduce the evolution parameter \cite{are/vol1} and consider the
dynamics on $E_{q}\left( 2\right) $ where (\ref{a}) are commutation
relations between momentum and coordinate, then we come to very complicated
and inappropriate conditions for the phase space. The reason lays in the
second term of the second formula in (\ref{a}). Such undesirable terms (so
called ``long'' solutions) result in difficulties while deriving and
exploiting of self-consistent Poisson bracket.

In search of ``short'' solutions we found the following ones on $E_q\left(
2\right) $ additionally to (\ref{a}) and (\ref{b}) 
\begin{equation}
\left\{ 
\begin{array}{lll}
\delta xx & = & x\delta x, \\ 
\delta xy & = & qy\delta x, \\ 
\delta yx & = & q^{-1}x\delta y, \\ 
\delta yy & = & y\delta y.
\end{array}
\right.  \label{c}
\end{equation}

In formulas (\ref{a}), (\ref{b}) and (\ref{c}) $\delta $ is the standard
exterior differential satisfying the classical Leibniz rule, lemma Poincare
and having the properties $\delta x\delta y=-q\delta y\delta x,\;\left(
\delta x\right) ^{2}=0,\;\left( \delta y\right) ^{2}=0$.

In contrast to the ``long'' solutions (\ref{a}) and (\ref{b}) which are
transformed one to another by the automorphism (\ref{auto}), the ``short''
solution (\ref{c}) is its fixed point.

\section{The $q$-deformed light-cone}

Let us consider $q$-deformed null-vector in 4-dimensional Minkowski space
which is described by the following $\left( 2\times 2\right) $ $q$-matrix 
\begin{equation}
X=\left( 
\begin{array}{cc}
x^{1\stackrel{.}{1}} & x^{1\stackrel{.}{2}} \\ 
x^{2\stackrel{.}{1}} & x^{2\stackrel{.}{2}}
\end{array}
\right) .  \label{x}
\end{equation}

We define the $q$-deformed light-cone as 
\begin{equation}
{\rm det}_{q^2}X=x^{1\stackrel{.}{1}}x^{2\stackrel{.}{2}}-q^2x^{1\stackrel{.%
}{2}}x^{2\stackrel{.}{1}}=0.  \label{det}
\end{equation}

Due to multiplicativity of the quantum determinant the Lorentz
transformations (represented by $SL_q\left( 2,C\right) $ matrices) do not
change the light-cone condition (\ref{det}).

The $q$-deformed components have the following commutation relations 
\begin{equation}
\left\{ 
\begin{array}{lll}
x^{1\stackrel{.}{1}}x^{1\stackrel{.}{2}} & = & q^2x^{1\stackrel{.}{2}}x^{1%
\stackrel{.}{1}}, \\ 
x^{1\stackrel{.}{1}}x^{2\stackrel{.}{1}} & = & q^2x^{2\stackrel{.}{1}}x^{1%
\stackrel{.}{1}}, \\ 
x^{1\stackrel{.}{2}}x^{2\stackrel{.}{1}} & = & x^{2\stackrel{.}{1}}x^{1%
\stackrel{.}{2}}, \\ 
x^{1\stackrel{.}{2}}x^{2\stackrel{.}{2}} & = & q^2x^{2\stackrel{.}{2}}x^{1%
\stackrel{.}{2}}, \\ 
x^{2\stackrel{.}{1}}x^{2\stackrel{.}{2}} & = & q^2x^{2\stackrel{.}{2}}x^{2%
\stackrel{.}{1}}, \\ 
x^{1\stackrel{.}{1}}x^{2\stackrel{.}{2}}-x^{2\stackrel{.}{2}}x^{1\stackrel{.%
}{1}} & = & \left( q^2-q^{-2}\right) x^{1\stackrel{.}{2}}x^{2\stackrel{.}{1}%
}.
\end{array}
\right.  \label{xx}
\end{equation}

As usual \cite{cartan,pen2} a null vector has a twistor representation. In $q
$-deformed case we introduce the twistors having $q$-deformed components 
\begin{equation}
\varphi _{q}=\left( \varphi _{q}^{A}\right) =\left( 
\begin{array}{c}
x \\ 
y
\end{array}
\right) ,\;A=1,2  \label{f}
\end{equation}

where $x$ and $y$ satisfy (\ref{q}).

The $q$-deformed Levi-Chivita tensor is 
\begin{equation}
\epsilon _q=\left( \epsilon _q^{AB}\right) =\left( 
\begin{array}{cc}
0 & q^{1/2} \\ 
-q^{-1/2} & 0
\end{array}
\right)  \label{e}
\end{equation}

and has the property 
\begin{equation}
\epsilon _q^{AB}\epsilon _{q,BA}=q+q^{-1}.  \label{eq}
\end{equation}

The $q$-antisymmetry of $\epsilon _q$ leads to 
\[
\varphi _q^A\varphi _q^B\epsilon _{q,AB}=\varphi _q^A\varphi _{q,A}=0. 
\]

The complex conjugated $q$-deformed twistor is 
\[
\overline{\varphi }_q=\left( \overline{\varphi }_q^{\stackrel{.}{A}}\right)
=\left( \overline{x},\overline{y}\right) . 
\]

We define here the involution $x\rightarrow \bar{x}$ as the standard
Hermitian conjugation.

The $q$-null-vector in terms of $q$-twistor components has the standard form 
\begin{equation}
X^{A\stackrel{.}{A}}=\varphi _q^A\overline{\varphi }_q^{\stackrel{.}{A}%
}=\left( 
\begin{array}{cc}
x\overline{x} & x\overline{y} \\ 
y\overline{x} & y\overline{y}
\end{array}
\right) .  \label{xxy}
\end{equation}

The light-cone condition (\ref{det}) in terms of $q$-twistor components is 
\begin{equation}
{\rm det}_{q^2}X=x\overline{x}y\overline{y}-q^2x\overline{y}y\overline{x}=0.
\label{xxyy}
\end{equation}

So we have two copies of the quantum plane $\left( x,y\right) $ and $\left( 
\bar{x},\bar{y}\right) $. Now we need such commutation relations between
components of the $q$-twistor and its conjugate which satisfy (\ref{xx}) and
the standard involution.

From (\ref{q}) and its conjugate we easily obtain the condition 
\begin{equation}
\bar{q}=q^{-1}=\exp \left( -ih\right)  \label{qe}
\end{equation}
or $\left| q\right| =1$.

Let other commutation relations have the general form 
\begin{equation}
\left\{ 
\begin{array}{ccc}
x\overline{x} & = & q^n\bar{x}x, \\ 
x\bar{y} & = & q^m\bar{y}x, \\ 
y\bar{x} & = & q^k\bar{x}y, \\ 
y\bar{y} & = & q^l\bar{y}y,
\end{array}
\right.  \label{xxqxx}
\end{equation}
where $n,m,k,l$ are arbitrary constants which should be determined from (\ref
{xx}) and involution. So we have the following independent equations 
\[
\left\{ 
\begin{array}{ccc}
m-n & = & 1, \\ 
n-k & = & 1, \\ 
l-k & = & 1.
\end{array}
\right. 
\]

Then 
\begin{equation}
\left\{ 
\begin{array}{ccc}
x\overline{x} & = & q^n\bar{x}x, \\ 
x\bar{y} & = & q^{n+1}\bar{y}x, \\ 
y\bar{x} & = & q^{n-1}\bar{x}y, \\ 
y\bar{y} & = & q^n\bar{y}y.
\end{array}
\right.  \label{xyqn}
\end{equation}

Now we consider the reality condition, as the consequence of the possibility
to reduce the dimension of $q$-deformed null-vector up to three where it is
real. So we derive $n=0$, and the commutation relations between twistor
components become 
\begin{equation}
\left\{ 
\begin{array}{lll}
xy & = & qyx, \\ 
x\overline{y} & = & q\overline{y}x, \\ 
\overline{x}y & = & qy\overline{x}, \\ 
\overline{x}\overline{y} & = & q\overline{y}\overline{x}, \\ 
\overline{x}x & = & x\overline{x}, \\ 
\overline{y}y & = & y\overline{y}.
\end{array}
\right.  \label{xy}
\end{equation}

Then we can find the ``short'' version of $q$-deformed differential calculus
on twistor components 
\begin{equation}
\left\{ 
\begin{array}{lll}
\delta xx & = & x\delta x, \\ 
\delta x\bar{x} & = & \bar{x}\delta x, \\ 
\delta xy & = & qy\delta x, \\ 
\delta x\bar{y} & = & q\bar{y}\delta x,
\end{array}
\right. \;\left\{ 
\begin{array}{lll}
\delta \bar{x}x & = & x\delta \bar{x}, \\ 
\delta \bar{x}\bar{x} & = & \bar{x}\delta \bar{x}, \\ 
\delta \bar{x}y & = & qy\delta \bar{x}, \\ 
\delta \bar{x}\bar{y} & = & q\bar{y}\delta \bar{x},
\end{array}
\right.  \label{dxy1}
\end{equation}
\begin{equation}
\left\{ 
\begin{array}{lll}
\delta yx & = & q^{-1}x\delta y, \\ 
\delta y\bar{x} & = & q^{-1}\bar{x}\delta y, \\ 
\delta yy & = & y\delta y, \\ 
\delta y\bar{y} & = & \bar{y}\delta y,
\end{array}
\right. \;\left\{ 
\begin{array}{lll}
\delta \bar{y}x & = & q^{-1}x\delta \bar{y}, \\ 
\delta \bar{y}\bar{x} & = & q^{-1}\bar{x}\delta \bar{y}, \\ 
\delta \bar{y}y & = & y\delta \bar{y}, \\ 
\delta \bar{y}\bar{y} & = & \bar{y}\delta \bar{y},
\end{array}
\right.  \label{dxy2}
\end{equation}
where $\left( \delta x\right) ^2=\left( \delta \bar{x}\right) ^2=\left(
\delta y\right) ^2=\left( \delta \bar{y}\right) ^2=0$.

The commutation relations between the differentials themselves take the form 
\begin{equation}
\left\{ 
\begin{array}{lll}
\delta x\delta \bar{x} & = & -\delta \bar{x}\delta x, \\ 
\delta x\delta y & = & -q\delta y\delta x, \\ 
\delta x\delta \bar{y} & = & -q\delta \bar{y}\delta x, \\ 
\delta y\delta x & = & -q^{-1}\delta x\delta y, \\ 
\delta y\delta \bar{x} & = & -q^{-1}\delta \bar{x}\delta y, \\ 
\delta y\delta \bar{y} & = & -\delta \bar{y}\delta y.
\end{array}
\right.  \label{dxdy}
\end{equation}

Using the obtained $q$-deformed differential calculus on $q$-twistors we can
build $q$-deformed differential calculus on any composite objects, as $q$%
-deformed null-vectors and tensors.

So for $q$-deformed null-vector we obtain 
\begin{equation}
\left\{ 
\begin{array}{lll}
\delta x^{1\stackrel{.}{1}}x^{1\stackrel{.}{1}} & = & x^{1\stackrel{.}{1}%
}\delta x^{1\stackrel{.}{1}}, \\ 
\delta x^{1\stackrel{.}{1}}x^{1\stackrel{.}{2}} & = & q^2x^{1\stackrel{.}{2}%
}\delta x^{1\stackrel{.}{1}}, \\ 
\delta x^{1\stackrel{.}{1}}x^{2\stackrel{.}{1}} & = & q^2x^{2\stackrel{.}{1}%
}\delta x^{1\stackrel{.}{1}}, \\ 
\delta x^{1\stackrel{.}{1}}x^{2\stackrel{.}{2}} & = & q^4x^{2\stackrel{.}{2}%
}\delta x^{1\stackrel{.}{1}},
\end{array}
\right. \;\left\{ 
\begin{array}{lll}
\delta x^{1\stackrel{.}{2}}x^{1\stackrel{.}{2}} & = & x^{1\stackrel{.}{2}%
}\delta x^{1\stackrel{.}{2}}, \\ 
\delta x^{1\stackrel{.}{2}}x^{2\stackrel{.}{1}} & = & x^{2\stackrel{.}{1}%
}\delta x^{1\stackrel{.}{2}}, \\ 
\delta x^{1\stackrel{.}{2}}x^{1\stackrel{.}{1}} & = & q^{-2}x^{1\stackrel{.}{%
1}}\delta x^{1\stackrel{.}{2}}, \\ 
\delta x^{1\stackrel{.}{2}}x^{2\stackrel{.}{2}} & = & q^2x^{2\stackrel{.}{2}%
}\delta x^{1\stackrel{.}{2}},
\end{array}
\right.  \label{dxx1}
\end{equation}
\begin{equation}
\left\{ 
\begin{array}{lll}
\delta x^{2\stackrel{.}{1}}x^{2\stackrel{.}{1}} & = & x^{2\stackrel{.}{1}%
}\delta x^{2\stackrel{.}{1}}, \\ 
\delta x^{2\stackrel{.}{1}}x^{1\stackrel{.}{1}} & = & q^{-2}x^{1\stackrel{.}{%
1}}\delta x^{2\stackrel{.}{1}}, \\ 
\delta x^{2\stackrel{.}{1}}x^{1\stackrel{.}{2}} & = & x^{1\stackrel{.}{2}%
}\delta x^{2\stackrel{.}{1}}, \\ 
\delta x^{2\stackrel{.}{1}}x^{2\stackrel{.}{2}} & = & q^2x^{2\stackrel{.}{2}%
}\delta x^{2\stackrel{.}{1}},
\end{array}
\right. \;\left\{ 
\begin{array}{lll}
\delta x^{2\stackrel{.}{2}}x^{2\stackrel{.}{2}} & = & x^{2\stackrel{.}{2}%
}\delta x^{2\stackrel{.}{2}}, \\ 
\delta x^{2\stackrel{.}{2}}x^{1\stackrel{.}{1}} & = & q^{-4}x^{1\stackrel{.}{%
1}}\delta x^{2\stackrel{.}{2}}, \\ 
\delta x^{2\stackrel{.}{2}}x^{1\stackrel{.}{2}} & = & q^{-2}x^{1\stackrel{.}{%
2}}\delta x^{2\stackrel{.}{2}}, \\ 
\delta x^{2\stackrel{.}{2}}x^{2\stackrel{.}{1}} & = & q^{-2}x^{2\stackrel{.}{%
1}}\delta x^{2\stackrel{.}{2}},
\end{array}
\right.  \label{dxx2}
\end{equation}

and 
\begin{equation}
\begin{array}{ccc}
\delta x^{1\stackrel{.}{1}}\delta x^{1\stackrel{.}{2}} & = & -q^2\delta x^{1%
\stackrel{.}{2}}\delta x^{1\stackrel{.}{1}}, \\ 
\delta x^{1\stackrel{.}{1}}\delta x^{2\stackrel{.}{1}} & = & -q^2\delta x^{2%
\stackrel{.}{1}}\delta x^{1\stackrel{.}{1}}, \\ 
\delta x^{1\stackrel{.}{1}}\delta x^{2\stackrel{.}{2}} & = & -q^4\delta x^{2%
\stackrel{.}{2}}\delta x^{1\stackrel{.}{1}}, \\ 
\delta x^{1\stackrel{.}{2}}\delta x^{2\stackrel{.}{1}} & = & -\delta x^{2%
\stackrel{.}{1}}\delta x^{1\stackrel{.}{2}}, \\ 
\delta x^{1\stackrel{.}{2}}\delta x^{2\stackrel{.}{2}} & = & -q^2\delta x^{2%
\stackrel{.}{2}}\delta x^{1\stackrel{.}{2}}, \\ 
\delta x^{2\stackrel{.}{1}}\delta x^{2\stackrel{.}{2}} & = & -q^2\delta x^{2%
\stackrel{.}{2}}\delta x^{2\stackrel{.}{1}}.
\end{array}
\label{qdxdy}
\end{equation}

From the above formulas we can find the commutation relations for
coordinates and derivatives 
\begin{equation}
\left\{ 
\begin{array}{ccc}
\frac \partial {\partial x^{1\stackrel{.}{1}}}x^{1\stackrel{.}{1}} & = & 
1+x^{1\stackrel{.}{1}}\frac \partial {\partial x^{1\stackrel{.}{1}}}, \\ 
\frac \partial {\partial x^{1\stackrel{.}{1}}}x^{1\stackrel{.}{2}} & = & 
q^{-2}x^{1\stackrel{.}{2}}\frac \partial {\partial x^{1\stackrel{.}{1}}}, \\ 
\frac \partial {\partial x^{1\stackrel{.}{1}}}x^{2\stackrel{.}{1}} & = & 
q^{-2}x^{2\stackrel{.}{1}}\frac \partial {\partial x^{1\stackrel{.}{1}}}, \\ 
\frac \partial {\partial x^{1\stackrel{.}{1}}}x^{2\stackrel{.}{2}} & = & 
q^{-4}x^{2\stackrel{.}{2}}\frac \partial {\partial x^{1\stackrel{.}{1}}},
\end{array}
\right. \;\left\{ 
\begin{array}{ccc}
\frac \partial {\partial x^{1\stackrel{.}{2}}}x^{1\stackrel{.}{1}} & = & 
q^2x^{1\stackrel{.}{1}}\frac \partial {\partial x^{1\stackrel{.}{2}}}, \\ 
\frac \partial {\partial x^{1\stackrel{.}{2}}}x^{1\stackrel{.}{2}} & = & 
1+x^{1\stackrel{.}{2}}\frac \partial {\partial x^{1\stackrel{.}{2}}}, \\ 
\frac \partial {\partial x^{1\stackrel{.}{2}}}x^{2\stackrel{.}{1}} & = & x^{2%
\stackrel{.}{1}}\frac \partial {\partial x^{1\stackrel{.}{2}}}, \\ 
\frac \partial {\partial x^{1\stackrel{.}{2}}}x^{2\stackrel{.}{1}} & = & 
q^{-2}x^{2\stackrel{.}{1}}\frac \partial {\partial x^{1\stackrel{.}{2}}},
\end{array}
\right.  \label{ddx1}
\end{equation}

\begin{equation}
\left\{ 
\begin{array}{ccc}
\frac \partial {\partial x^{2\stackrel{.}{1}}}x^{1\stackrel{.}{1}} & = & 
q^2x^{1\stackrel{.}{1}}\frac \partial {\partial x^{2\stackrel{.}{1}}}, \\ 
\frac \partial {\partial x^{2\stackrel{.}{1}}}x^{2\stackrel{.}{1}} & = & 
1+x^{2\stackrel{.}{1}}\frac \partial {\partial x^{2\stackrel{.}{1}}}, \\ 
\frac \partial {\partial x^{2\stackrel{.}{1}}}x^{1\stackrel{.}{2}} & = & x^{1%
\stackrel{.}{2}}\frac \partial {\partial x^{2\stackrel{.}{1}}}, \\ 
\frac \partial {\partial x^{2\stackrel{.}{1}}}x^{2\stackrel{.}{2}} & = & 
q^{-2}x^{2\stackrel{.}{2}}\frac \partial {\partial x^{2\stackrel{.}{1}}},
\end{array}
\right. \;\left\{ 
\begin{array}{ccc}
\frac \partial {\partial x^{2\stackrel{.}{2}}}x^{1\stackrel{.}{1}} & = & 
q^4x^{1\stackrel{.}{1}}\frac \partial {\partial x^{2\stackrel{.}{2}}}, \\ 
\frac \partial {\partial x^{2\stackrel{.}{2}}}x^{1\stackrel{.}{2}} & = & 
q^2x^{1\stackrel{.}{2}}\frac \partial {\partial x^{2\stackrel{.}{2}}}, \\ 
\frac \partial {\partial x^{2\stackrel{.}{2}}}x^{2\stackrel{.}{1}} & = & 
q^2x^{2\stackrel{.}{1}}\frac \partial {\partial x^{2\stackrel{.}{2}}}, \\ 
\frac \partial {\partial x^{2\stackrel{.}{2}}}x^{2\stackrel{.}{2}} & = & 
1+x^{2\stackrel{.}{2}}\frac \partial {\partial x^{2\stackrel{.}{2}}}
\end{array}
\right.  \label{ddx2}
\end{equation}

and between the derivatives 
\begin{equation}
\begin{array}{ccc}
\frac \partial {\partial x^{1\stackrel{.}{1}}}\frac \partial {\partial x^{1%
\stackrel{.}{2}}} & = & q^2\frac \partial {\partial x^{1\stackrel{.}{2}%
}}\frac \partial {\partial x^{1\stackrel{.}{1}}}, \\ 
&  &  \\ 
\frac \partial {\partial x^{1\stackrel{.}{1}}}\frac \partial {\partial x^{2%
\stackrel{.}{1}}} & = & q^2\frac \partial {\partial x^{2\stackrel{.}{1}%
}}\frac \partial {\partial x^{1\stackrel{.}{1}}}, \\ 
&  &  \\ 
\frac \partial {\partial x^{1\stackrel{.}{1}}}\frac \partial {\partial x^{2%
\stackrel{.}{2}}} & = & q^4\frac \partial {\partial x^{2\stackrel{.}{2}%
}}\frac \partial {\partial x^{1\stackrel{.}{1}}}, \\ 
&  &  \\ 
\frac \partial {\partial x^{1\stackrel{.}{2}}}\frac \partial {\partial x^{2%
\stackrel{.}{1}}} & = & \frac \partial {\partial x^{2\stackrel{.}{1}}}\frac
\partial {\partial x^{1\stackrel{.}{2}}}, \\ 
&  &  \\ 
\frac \partial {\partial x^{1\stackrel{.}{2}}}\frac \partial {\partial x^{2%
\stackrel{.}{2}}} & = & q^2\frac \partial {\partial x^{2\stackrel{.}{2}%
}}\frac \partial {\partial x^{1\stackrel{.}{2}}}, \\ 
&  &  \\ 
\frac \partial {\partial x^{2\stackrel{.}{1}}}\frac \partial {\partial x^{2%
\stackrel{.}{2}}} & = & q^2\frac \partial {\partial x^{2\stackrel{.}{2}%
}}\frac \partial {\partial x^{2\stackrel{.}{1}}}.
\end{array}
\label{qqdxdx}
\end{equation}

By analogy we can write the commutation relation between $q$-deformed
differentials and derivatives.

The above relations are consistent with the following condition 
\begin{equation}
\delta \left( {\rm det}_{q^2}X\right) =0.  \label{ddet}
\end{equation}

\section{Momenta and $q$-D'Alembertian}

Let us introduce $q$-matrix for momenta 
\begin{equation}
P_q=\left( 
\begin{array}{cc}
P_{1\stackrel{.}{1}} & P_{1\stackrel{.}{2}} \\ 
P_{2\stackrel{.}{1}} & P_{2\stackrel{.}{2}}
\end{array}
\right) =\left( 
\begin{array}{cc}
-i\frac \partial {\partial x^{1\stackrel{.}{1}}} & -i\frac \partial
{\partial x^{1\stackrel{.}{2}}} \\ 
&  \\ 
-i\frac \partial {\partial x^{2\stackrel{.}{1}}} & -i\frac \partial
{\partial x^{2\stackrel{.}{2}}}
\end{array}
\right)  \label{pq}
\end{equation}
with commutation relations determined by (\ref{qqdxdx}).

The $q$-D'Alembertian is defined by

\begin{equation}
\Box _{q^2}=-{\rm det}_{q^2}P_q=\frac \partial {\partial x^{1\stackrel{.}{1}%
}}\frac \partial {\partial x^{2\stackrel{.}{2}}}-q^2\frac \partial {\partial
x^{1\stackrel{.}{2}}}\frac \partial {\partial x^{2\stackrel{.}{1}}}.
\label{pdet}
\end{equation}

The light-cone condition (\ref{det}) leads to the similar condition for
momenta 
\begin{equation}
\Box _{q^2}=0.  \label{dal0}
\end{equation}

In the classical limit $q=\exp \left( ih\right) \rightarrow 1$ we decompose
the $q$-D'Alembertian as follows 
\begin{eqnarray}
\Box _{q^{2}} &=&\Box -\left( q^{2}-1\right) \frac{\partial }{\partial x^{1%
\stackrel{.}{2}}}\frac{\partial }{\partial x^{2\stackrel{.}{1}}}  \nonumber
\\
&=&\Box -2ihP_{1\stackrel{.}{2}}P_{2\stackrel{.}{1}}=0,  \label{dalclas}
\end{eqnarray}
where $\Box $ is the ordinary D'Alembertian.

The second term in (\ref{dalclas}) gives us a new way of the off-shell
approximation and is responsible for its exact shape. In contrast to the
standard picture \cite{azc/kul/rod} the momenta entering into the additional
off-shell term in (\ref{dalclas}) commute.

\section{Conclusion}

We have proposed a version of the differential calculus on $q$-deformed
light-cone which can be applied to description of the dynamics of the
massless quantum particles, $q$-deformed null-strings and null-membranes.

Using the obtained $q$-deformed differential calculus on $q$-twistors we
have the possibility to construct the corresponding calculi on $q$-tensors
of any rank.

\section*{Acknowledgments}

V. P. Akulov is grateful to A. P. Isaev, P. N. Pyatov and B. M. Zupnik for
discussions.

\end{document}